\documentclass{article}

\usepackage{microtype}
\usepackage{graphicx}
\usepackage{booktabs} 
\usepackage{bm}
\usepackage{amsthm}
\usepackage{lipsum,multicol}
\usepackage{amsfonts}
\usepackage[mathscr]{euscript}
\usepackage{url}
\usepackage{eepic, epsfig, amsmath, amssymb, amsthm,latexsym, setspace}
\usepackage{rotating}
\usepackage{nicefrac}
\usepackage{lipsum,multicol}
\usepackage[parfill]{parskip}
\usepackage{caption}
\usepackage{subcaption}
\usepackage{float}
\usepackage{algorithm}
\usepackage{algorithmic}
\usepackage{times}
\usepackage{xcolor}
\usepackage{cite}
\usepackage{geometry}
\usepackage[shortcuts]{extdash}

\numberwithin{equation}{section}

\newcommand{\prox}{{\rm{prox}}}

\newcommand{\mcG}{\mathcal{G}}
\newcommand{\mcS}{\mathcal{S}}

\newcommand{\mbR}{\mathbb{R}}

\usepackage{hyperref}

\begin{document}
	
	\title{Solving Large Scale Quadratic Constrained Basis Pursuit}
	\date{\today}
	\author{Jirong Yi\thanks{Department of Electrical and Computer Engineering, University of Iowa, Iowa City, USA. All copyrights are reserved.}
	}

	\maketitle
	
%
%
%

We solve
\begin{align}\label{Defn:QCBP}
	\min_x \|x\|_1, {\rm s.t.\ }\|y-Ax\|_2\leq\eta,
\end{align}
where $y\in\mbR^m, A\in\mbR^{m\times d} (m<d)$, and $\eta>0$ are known. We will adopt the alternating direction method of multipliers (ADMM) and the idea of operator splitting to design efficient algorithm for solving the above quadratically constrained basis pursuit problem \cite{boyd_distributed_2010,fougner_parameter_2015,parikh_block_2014,parikh_proximal_2013}.

\section{Theoretical guarantees}

We reformulate \eqref{Defn:QCBP} as
\begin{align}\label{Defn:QCBPinGraph}
	\min_{x,z} g(x) + f(z), {\rm \ s.t.\ }Ax=z,
\end{align}
where $g(x) = \|x\|_1$, and $f(z):\mbR^m \to \mbR$ is an indicator function defined as
\begin{align}\label{Defn:fzIndicator}
	f(z):=
	\begin{cases}
	0, z\in\Omega,\\
	\infty, z\notin\Omega,
	\end{cases}
	\Omega = \{z: \|z-y\|_2\leq \eta\}.
\end{align}

The Lagrangian dual of \eqref{Defn:fzIndicator} is
\begin{align}\label{Defn:LagrangianDual}
	\mathcal{L}(v)
	&:= \min_{x,z} \left(g(x) + f(z) + v^T(Ax-z)\right)\nonumber\\
	& = - f^*(v) - g^*(-A^Tv),
\end{align}
where $v\in\mbR^m$ is the dual variable, $f^*(v)$ is the convex conjugate of $f(z)$ at $v\in\mbR^m$, and $g^*(-A^Tv)$ is the convex conjugate of $g(x)$ at $-A^Tv$, i.e., 
\begin{align}\label{Eq:gx_Conjugate}
	g^*(-A^Tv): = 
	\begin{cases}
	0, -A^Tv \in \Omega',\\
	\infty, -A^Tv \notin \Omega',\\
	\end{cases}
	\Omega':=\{x: \|x\|_\infty \leq 1\}
\end{align}
and
\begin{align}\label{Eq:fz_Conjugate}
	f^*(v):= y + \frac{\eta}{\|v\|_2} v.
\end{align}
The dual problem can be formulated as
\begin{align}
	\max_{v,\mu} -f^*(v) - g^*(\mu), {\rm\ s.t.\ } -A^Tv = \mu.
\end{align}

Assume the Slater's condition holds, i.e., there exists $x\in\mbR^d$ such that $\|y-Ax\|_2<\eta$, then the convexity of problem \eqref{Defn:QCBP} implies that the optimal solution will achieve zero duality gap, i.e.,
\begin{align}\label{Eq:ZeroDuality}
	g(x) + f(z) = - f^*(v) - g^*(\mu).
\end{align}
From KKT conditions, the optimal solution must satisfy \eqref{Eq:ZeroDuality} and 
\begin{align}\label{Eq:Feasibility}
	Ax = z, -A^Tv = \mu.
\end{align}
Thus, the \eqref{Eq:ZeroDuality} and \eqref{Eq:Feasibility} can be used as optimality certificates or stopping criterion in algorithm design. More specifically, we define primal residual, dual residual, and duality gap with respect to a certain tuple $(x,z,v,\mu)$ as
\begin{align}\label{Defn:StoppingCriterion}
	r_p := \|Ax-z\|,
	r_d:=\|A^Tv + \mu\|,
	\delta_g = g(x) + f(z) + f^*(v) + g^*(\mu).
\end{align}

\section{Algorithm design based on ADMM}

We adopt ideas from alternating projection methods, and reformulate \eqref{Defn:QCBPinGraph} as 
\begin{align}\label{Defn:QCBPinADMM}
	\min_{x,z,x',z'} g(x) + f(z) + I_\mathcal{G}(x',z'), {\rm\ s.t.\ } 
	\left[\begin{matrix}
	x \\ z
	\end{matrix}\right]
	=
	\left[\begin{matrix}
	x' \\ z'
	\end{matrix}\right],
\end{align}
where $I_\mathcal{G}(x',z')$ is defined as
\begin{align}\label{Defn:GraphSet}
	I_\mathcal{G}(x',z'):=
	\begin{cases}
	0, (x',z') \in\mcG,\\
	\infty, (x',z') \notin \mcG,
	\end{cases}
	\mcG:=\{(x,z): Ax = z\}.
\end{align}
The augmented Lagrangian of \eqref{Defn:GraphSet} becomes
\begin{align}\label{Defn:AugmentedLagrangian}
	L_\rho(x,y,x',y',v) 
	:= g(x) + f(z) + I_\mathcal{G}(x',z') 
	+ V^T\left(\left[\begin{matrix}
	x \\ z
	\end{matrix}\right]
	-
	\left[\begin{matrix}
	x' \\ z'
	\end{matrix}\right]\right)
	+ \frac{\rho}{2} \left\|\left[\begin{matrix}
	x \\ z
	\end{matrix}\right]
	-
	\left[\begin{matrix}
	x' \\ z'
	\end{matrix}\right]\right\|_2^2,
\end{align}
where $V:=[v_x^T\ v_z^T]^T\in\mbR^{d+m}$ is the dual variable, and $\rho>0$ is a parameter. Define
\begin{align}\label{Defn:GroupedVars}
X:=\left[\begin{matrix}
x \\ z
\end{matrix}\right] \in\mbR^{d+m},
X':=\left[\begin{matrix}
x' \\ z'
\end{matrix}\right] \in\mbR^{d+m},
\end{align}
and we get the iterations in ADMM are
\begin{align}\label{Defn:ADMMiterates}
\begin{cases}
X^{k+1}:= \arg\min_{X} L_\rho(X,{X'}^{k},V^{k}),\\
{X'}^{k+1}:= \arg\min_{X'} L_\rho(X^{k+1},{X'},V^{k}),\\
v^{k+1}:=V^k + \rho\left(X^{k+1} - {X'}^{k+1}\right).
\end{cases} 
\end{align} 

More specifically, 
\begin{align}\label{Defn:ADMMiterates1}
X^{k+1}
&:= \arg\min_{X} L_\rho(X,{X'}^{k},V^{k}) \nonumber\\
& = \left[\begin{matrix}
\arg\min_x \left(g(x) + {v_x^{k}}^T(x-{x'}^{k})  + \frac{\rho}{2} \|x-{x'}^{k}\|_2^2 \right) \\
\arg\min_z \left(f(z) + {v_z^{k}}^T(z - {z'}^{k}) + \frac{\rho}{2} \|z - {z'}^{k}\|_2^2\right)
\end{matrix}\right] \nonumber\\
& = \left[\begin{matrix}
\arg\min_x \left(g(x) + \frac{\rho}{2} \left\|x - \left({x'}^{k} - \frac{v_x^{k}}{\rho}\right)\right\|_2^2 \right) \\
\arg\min_z \left(f(z) + \frac{\rho}{2} \left\|z - \left( {z'}^{k} - \frac{v_z^{k}}{\rho}\right)\right\|_2^2\right)
\end{matrix}\right]
\end{align} 
or simply
\begin{align}\label{Eq:Defn:ADMMiterates1inProximatorForm}
	X^{k+1} = \left[\begin{matrix}
	\prox_g({x'}^{k} - \tilde{x}^{k}) \\
	\prox_f({z'}^{k} - \tilde{z}^{k})
	\end{matrix}\right]
\end{align}
where $\prox_g(v)$ is the proximator of function $g$ at $v$ which is defined as
\begin{align}
     \prox_g(v) 
     & = \arg\min_x \left(g(x) + \frac{\rho}{2} \left\|x - v\right\|_2^2 \right) 
\end{align}
The $\tilde{x}^{k}$ and $\tilde{z}^{k}$ are defined as
\begin{align}\label{Defn:ScaledDualVariable}
	\tilde{x}^{k}:=\frac{v_x^{k}}{\rho}, 
	\tilde{z}^{k}:=\frac{v_z^{k}}{\rho}.
\end{align}
More specifically, the proximator of $g$ at $v\in\mbR^d$ is 
\begin{align}\label{Eq:ProximatorL1}
	\prox_g(v) = \mcS_{1/\rho}(v),
\end{align}
where $\mcS_{1/\rho}(\cdot)$ is the elementwise soft thresholding function, i.e., 
\begin{align}
	\left[\mcS_{1/\rho}(v)\right]_i :=
	\begin{cases}
	v_i - 1/\rho, v_i>1/\rho,\\
	0, |v_i| \leq 1/\rho,\\
	v_i + 1/\rho, v_i < -1/\rho.
	\end{cases} 
\end{align}
The proximator of $f$ at $v\in\mbR^m$ is
\begin{align}\label{Eq:ProximatorBallIndicator}
	\prox_f(v) = \frac{\eta}{\|v-y\|_2} (v-y) + y.
\end{align}

The updating rule for ${X'}$ can be specified as
\begin{align}\label{Defn:ADMMiterates2}
	{X'}^{k+1}
	& :=\arg\min_{X'} L_\rho(X^{k+1},{X'},V^{k}) \nonumber\\
	& = \arg\min_{X' \in \mcG} v^T\left(\left[\begin{matrix}
	x \\ z
	\end{matrix}\right]
	-
	\left[\begin{matrix}
	x' \\ z'
	\end{matrix}\right]\right)
	+ \frac{\rho}{2} \left\|\left[\begin{matrix}
	x \\ z
	\end{matrix}\right]
	-
	\left[\begin{matrix}
	x' \\ z'
	\end{matrix}\right]\right\|_2^2 \nonumber\\
	& = \arg\min_{X'\in\mcG} \frac{1}{2}\left\|
	X' - \left(X^{k+1} + \frac{v^{k}}{\rho}\right)
	\right\|_2^2 \nonumber\\
	& = \arg\min_{X'\in\mcG} \frac{1}{2}\left\|x' - (x^{k+1} + \tilde{x}^{k})\right\|_2^2 
	+ \frac{1}{2}\left\|z' - (z^{k+1} + \tilde{z}^{k})\right\|_2^2 \nonumber\\
	& = \prod(x^{k+1} + \tilde{x}^{k},z^{k+1} + \tilde{z}^{k})
\end{align}
where $\prod(x,z)$ is the projection of $(x,z)$ onto $G$, i.e., the solution to
\begin{align}\label{Defn:ADMMiterates2inProjectionForm}
	\min_{x',z'} \frac{1}{2}\|x'-x\|_2^2 + \frac{1}{2}\|z'-z\|_2^2, {\rm\ s.t.\ }Ax'=z'.
\end{align}

Define
\begin{align}\label{Defn:ScaledDual}
	\tilde{X}:=\frac{V}{\rho} = \left[\begin{matrix}
	\frac{v_x}{\rho} \\ \frac{v_z}{\rho}
	\end{matrix}\right],
\end{align}
and the updating rule for dual variable $V$ can be written as
\begin{align}\label{Eq:ADMMiterates3inScaled}
	\tilde{X}^{k+1} = \tilde{X}^{k} + \left(X^{k+1} - {X'}^{k}\right).
\end{align}

\subsection{Analytic solution to \eqref{Defn:ADMMiterates2inProjectionForm}}

Since \eqref{Defn:ADMMiterates2inProjectionForm} is convex, from the KKT conditions of \eqref{Defn:ADMMiterates2inProjectionForm}, we know that $x',z'$ are the optimal solution to \eqref{Defn:ADMMiterates2inProjectionForm} if and only if there exists $\lambda\in\mbR^m$ such that
\begin{align}
    \begin{cases}
    Ax' = z',\\
    x'-c + A^T\lambda= 0,\\
    z'-d-\lambda = 0,
    \end{cases}
\end{align}
which implies that the optimal $x',z'$ can be obtained from solving the following linear system
\begin{align}
	\left[\begin{matrix}
	A & -I_{m\times m} \\
	I_{d\times d} & A^T
	\end{matrix}\right]
	\left[\begin{matrix}
	x' \\ z'
	\end{matrix}\right]
	 = \left[\begin{matrix}
	 0 \\ A^Tz + x
	 \end{matrix}\right]
\end{align}

{\bf Remarks:} (1) the matrix $\left[\begin{matrix}
A & -I_{m\times m} \\
I_{d\times d} & A^T
\end{matrix}\right]$ is highly sparse, and this structure can be combined with other potential structured of $A$ to simplify the computation; (2) even simple elimination can be used to simplify the problem, i.e.,
\begin{align}\label{Eq:EliminationForm1}
	\begin{cases}
	z' = (AA^T+I_{m\times m})^{-1}\left(AA^Tz + Ax\right),\\
	x' = A^T(z-z') + x,
	\end{cases}
\end{align}
or
\begin{align}\label{Eq:EliminationForm2}
	\begin{cases}
	x'=(A^TA+I_{d\times d})^{-1}(x+A^Tz),\\
	z'=Ax'.
	\end{cases}
\end{align}
Both the two matrices $AA^T+I_{m\times m}$ and $A^TA+I_{d\times d}$ are positive definite, thus factorization techniques can be used to accelerate the computation; (3) since $m<d$, the \eqref{Eq:EliminationForm1} will be more efficient; (4) apply Cholesky decomposition once to get $AA^T+I_{m\times m}=LL^T$; (5) calculate $AA^T$ once; (6) solve for $z'$ backward, i.e., $z'=L^{-T}L^{-1}(AA^Tz+Ax)$;  

\subsection{Algorithm in pseudocodes}\label{Sec:Pseudocodes}

The algorithm can be summarized as in Algorithm \ref{Alg:QCBPviaADMM}.

Computational complexity - running time: (1) line 5, 7, and 8 takes $O(d+m)$; (2) line 6 takes $O(dm^2)$ for Cholesky decomposition over $AA^T+I_{m\times m}$, $O(dm^2)$ for $AA^T$ once, $O(dm)$ for backward solving $z'$ using \eqref{Eq:EliminationForm1}; (3) line 9 and 10 takes $O(dm)$. Thus, $O(dm^2)$ but only once in total; 

Computational complexity - space or memory: $O(m^2)$; 

Baseline algorithm, CVX using interior point method: (1) $O(md^2)$ but multiple times.    

\begin{algorithm}[!htb]
	\caption{Algorihm for solving large scale QCBP}\label{Alg:QCBPviaADMM}
	\begin{algorithmic}[1]
		
		\STATE {\bf Input:} $A\in\mbR^{m\times d}, y\in\mbR^m$, and $\eta>0$
		
		\STATE {\bf Parameters:} $\rho>0$, $MaxIte\in\mathbb{Z}_+$, $\epsilon_p>0$, $\epsilon_d>0$, and $\epsilon_{dg}>0$
		
		\STATE {\bf Initialization:} $k=0$, $X^{0}={X'}^{0}=\left[\begin{matrix}
		{x'}^{0} \\ {z'}^{0}
		\end{matrix}\right]\in\mbR^{d+m}$ 
		and 
		$\tilde{X}^{0} = \left[\begin{matrix}
		\tilde{x}^{0} \\ \tilde{z}^{0}
		\end{matrix}\right]\in\mbR^{d+m}$
		
		\WHILE{$k\leq MaxIte$}
		
		\STATE Solve $X^{k+1}$ via \eqref{Eq:Defn:ADMMiterates1inProximatorForm}, i.e., $	X^{k+1} = \left[\begin{matrix}
		\prox_g({x'}^{k} - \tilde{x}^{k}) \\
		\prox_f({z'}^{k} - \tilde{z}^{k})
		\end{matrix}\right]$
		
		\STATE Solve ${X'}^{k+1}$ via \eqref{Defn:ADMMiterates2}, i.e., ${X'}^{k}=\prod(x^{k+1} + \tilde{x}^{k},z^{k+1} + \tilde{z}^{k})$ via \eqref{Eq:EliminationForm1}
		
		\STATE Solve $\tilde{X}^{k+1}$ via \eqref{Eq:ADMMiterates3inScaled}, i.e., $\tilde{X}^{k+1} = \tilde{X}^{k} + \frac{1}{\rho}\left(X^{k+1} - {X'}^{k}\right)$
		
		\STATE Get $V^{k+1}$ via \eqref{Defn:ScaledDual}, i.e., $V^{k+1} = \rho \tilde{X}^{k+1}$
		
		\STATE Calculate primal residual via \eqref{Defn:StoppingCriterion}, i.e., $r_p^{k+1} = \|Ax^{k+1} - z^{k+1}\|_2$
		
		\STATE Calculate dual residual via \eqref{Defn:StoppingCriterion}, i.e., $r_d^{k+1} = \|A^Tv_z^{k+1} + v_x^{k+1}\|_2$
		
		\STATE Calculate duality gap via \eqref{Defn:StoppingCriterion}, i.e., $\delta_g^{k+1} = g(x^{k+1}) + f(z^{k+1}) + f^*(v_z^{k+1}) + g^*(v_x^{k+1})$
		
		\IF{$r_p^{k+1}<=\epsilon_p$ and $r_d^{k+1}<=\epsilon_d$ and $\delta_g^{k+1}\leq \epsilon_{dg}$}
		
		    \STATE break
		    
		\ELSE
		
		    \STATE $k=k+1$
		\ENDIF
		
		\ENDWHILE
		
		\IF{$k\geq MaxIter$}
		
		    \STATE Algorithm does not converge in $MaxIter$ iterations
		    
		    \STATE {\bf Return} NOT CONVERGED
		
		\ELSE
		    
		    \STATE Algorithm converges within $MaxIter$ iterations
		    
		    \STATE {\bf Return} $X^{k+1}$
		
		\ENDIF
		
	\end{algorithmic}
\end{algorithm}

\section{Numerical experiments}

Computational environment: (1) desktop with Intel(R) Core(TM) i7-6700 CPU \@ 3.40GHz   3.40 GHz, 32.0 GB RAM; (2) OS Windows 10 Education; (3) MATLAB R2018a; (4) baseline CVX which solves \eqref{Defn:QCBP} using interior point method; 

Computational setup: (1) $x$ is assumed to be sparse with cardinality $k=p_s*d$, and generated randomly; (2) $m=p_m*d$, and generate $A$ randomly; (3) generate noise $v\in\mbR^m$ randomly and normalize it to have magnitude $\eta$; (4) $y$ is assumed to be generated via $y=Ax+v$; 

\begin{figure}
	\centering
	\includegraphics[width=\linewidth]{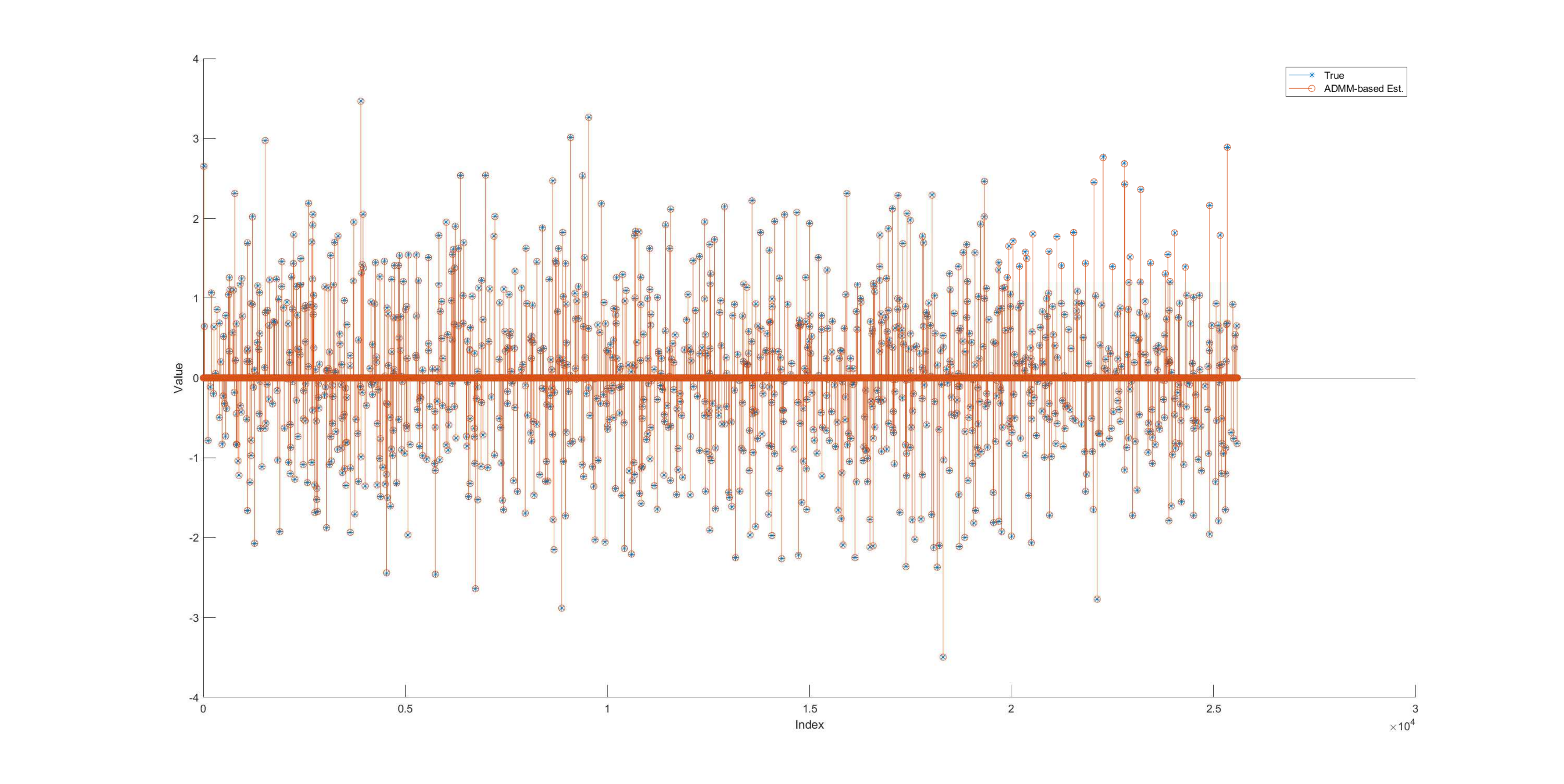}
	\caption{$p_s=0.4$, $p_m=0.05$, $\eta=0.1$, $d=25600$}\label{Fig:DiffSize}
\end{figure}

Results: see Table \ref{Tab:DiffSize} and Figure \ref{Fig:DiffSize}

\begin{table}[!htb]
	\centering
	\begin{tabular}{|l|l|l|l|l|l|}
		\hline
	Time (sec)	& $d=$100 & $d=$400 & $d=$1600 & $d=$6400 & $d=$25600 \\
	\hline
	CVX	& 0.7 & 1.16 & 44 & NA & NA \\
	\hline
	Algorithm \ref{Alg:QCBPviaADMM}	& 0.01 & 0.02 & 0.31 & 4.82 & 104.91\\ 
	\hline
	\end{tabular}
\caption{Computational performance comparisons: $p_s=0.4$, $p_m=0.05$, $\eta=0.1$}\label{Tab:DiffSize}
\end{table}

\bibliographystyle{unsrt}
\bibliography{Ref_qcbp}

\begin{thebibliography}{1}

\bibitem{boyd_distributed_2010}
S.~Boyd, N.~Parikh, E.~Chu, B.~Peleato, and J.~Eckstein.
\newblock Distributed optimization and statistical learning via the alternating
  direction method of multipliers.
\newblock {\em Machine Learning}, 3(1):1--122, 2010.

\bibitem{fougner_parameter_2015}
C.~Fougner and S.~Boyd.
\newblock Parameter selection and pre-conditioning for a graph form solver.
\newblock {\em arXiv:1503.08366 [math]}, March 2015.
\newblock arXiv: 1503.08366.

\bibitem{parikh_block_2014}
Neal Parikh and Stephen Boyd.
\newblock Block splitting for distributed optimization.
\newblock {\em Mathematical Programming Computation}, 6(1):77--102, 2014.
\newblock Publisher: Springer.

\bibitem{parikh_proximal_2013}
N.~Parikh and S.~Boyd.
\newblock Proximal algorithms.
\newblock {\em Foundations and Trends in Optimization}, 1(3):123--231, 2013.

\end{thebibliography}

\end{document}